\documentclass[fleqn,10pt]{wlscirep}
\usepackage[utf8]{inputenc}
\usepackage{graphicx}
\usepackage[T1]{fontenc}
 \usepackage[normalem]{ulem} 
\usepackage{amssymb,amsfonts,amsmath}

\usepackage{url}

\begin{document}

\title{The Waiting-Time Paradox}
\author[1,2,*]{Naoki Masuda}
\author[3]{Mason A. Porter}
\affil[1]{Department of Mathematics,
State University of New York at Buffalo, Buffalo, NY, United States}
\affil[2]{Computational and Data-Enabled Science and Engineering Program, State University of New York at Buffalo, Buffalo, NY, United States}
\affil[3]{Department of Mathematics, University of California Los Angeles, Los Angeles, CA, United States}
\affil[*]{Correspondence to naokimas@buffalo.edu 
}


\begin{abstract}
Suppose that you're going to school and arrive at a bus stop. How long do you have to wait before the next bus arrives? Surprisingly, it is longer --- possibly much longer --- than what the bus schedule suggests intuitively. This phenomenon, which is called the \emph{waiting-time paradox}, has a purely mathematical origin. Different buses arrive with different intervals, leading to this paradox. In this article, we explore the waiting-time paradox, explain why it happens, and discuss some of its implications (beyond the possibility of being late for school).
\end{abstract}

\maketitle





\section{How long do you have to wait for the next bus?}

\begin{figure}[h]
\includegraphics[width=16cm]{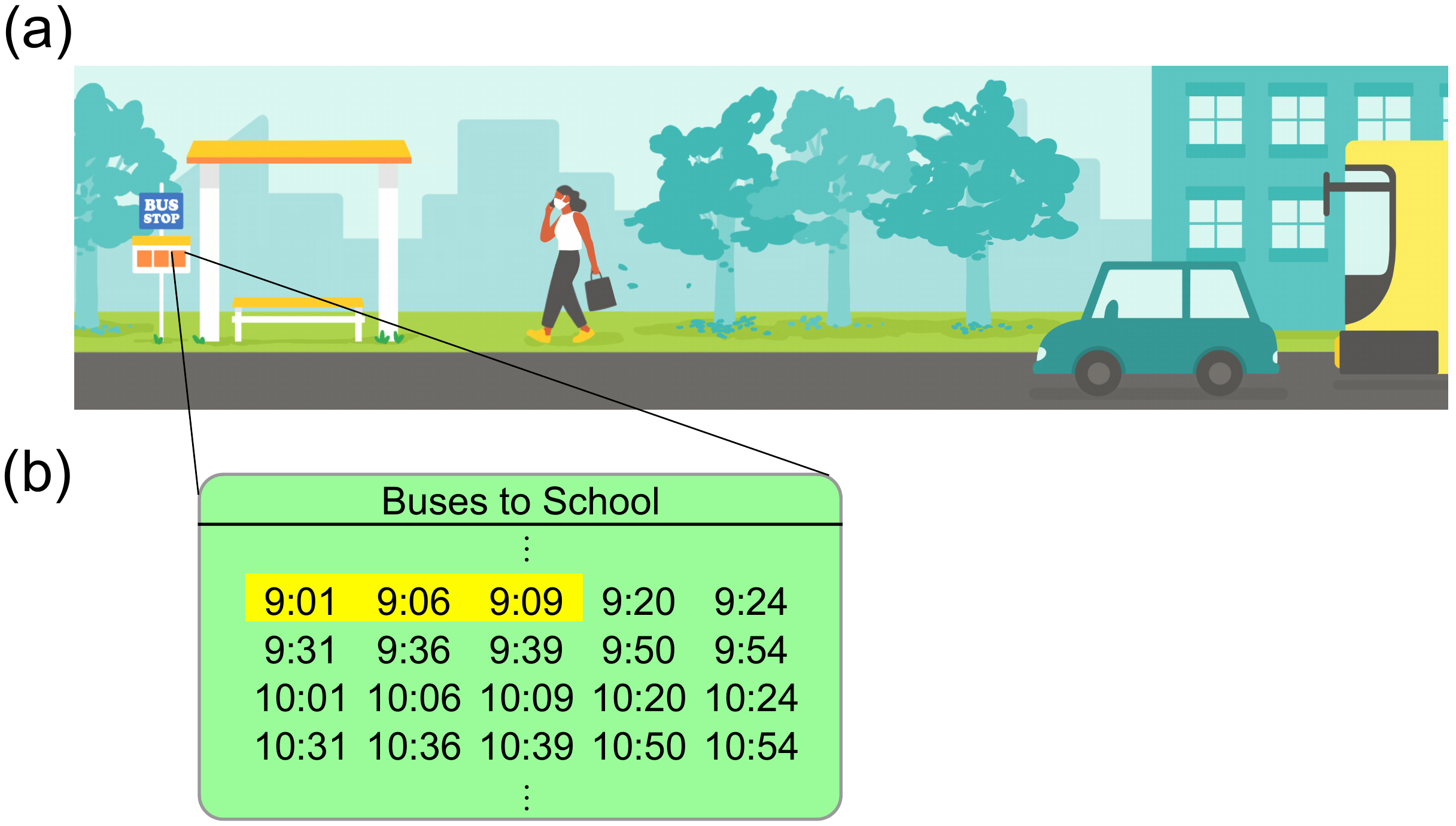}
\caption{
(a) You go to a bus stop and wait for the next bus. (b) You look at the schedule and see that there are 10 buses each hour. Therefore, there is 1 bus every 6 minutes on average.
The first three buses after 9:00 are scheduled at 9:01, 9:06, and 9:09, which we highlight in yellow in (b). The associated inter-event times --- that is, the times between consecutive buses --- are 5 minutes (between the first and second buses) and 3 minutes (between the second and third buses). 
[Panel (a) was drawn by Iris Leung.]
}
\label{fig:bus example}
\end{figure}

Suppose that you live in a city and take a bus to go to school. Because buses come frequently in your neighborhood, perhaps you do not need to pay close attention to the bus schedule. Maybe you just go to your bus stop and ride the next bus [see Figure~\ref{fig:bus example}(a)]. If you arrive and have no idea when the next bus is coming (perhaps you forgot to look at the schedule), how long do you have to wait for the next bus?

Now suppose that 10 buses come each hour [see Figure~\ref{fig:bus example}(b)], so 1 bus comes every 6 minutes on average. If the most recent bus left right before you arrive, you may have to wait a while for the next bus (perhaps 5 or 6 minutes). Alternatively, if the next bus left earlier than that, maybe the next bus will arrive in only 1 minute? Otherwise, the next bus may arrive in 4 minutes, or perhaps 2 minutes?
An educated guess for how long you have to wait --- your ``waiting time'' --- is 3 minutes, which is half the time between buses on average.

%

However, the above reasoning is not correct. Typically, you actually have to wait much longer than 3 minutes. Your expected waiting time can be even longer than 6 minutes! This phenomenon is called the \textbf{waiting-time paradox}, which is also called the ``bus paradox'' and the ``inspection paradox'' \cite{Weldman1957JOperResSoc,Masuda2020NortheastJCompSyst}. Here is why people think of it as a paradox: your typical waiting time at the bus station is longer than half of the average interval of time between buses (which is 3 minutes in the example above). It is a mathematical phenomenon and has nothing to do with buses or schedules or the navigation strategy of a bus driver. In university courses and scientific research, the waiting-time paradox shows up a lot in topics like \textbf{probability theory}, queuing theory, and network analysis. We will briefly discuss some examples later in the paper.


\section{The mathematics behind the waiting-time paradox}

\subsection{Simplifying a bus schedule to understand it more easily}

Why does the bus paradox occur? A complete explanation calls for some probability theory, but one can understand the key logic behind the paradox without advanced mathematics. Let's give it a try.

\begin{figure}
\includegraphics[width=16cm]{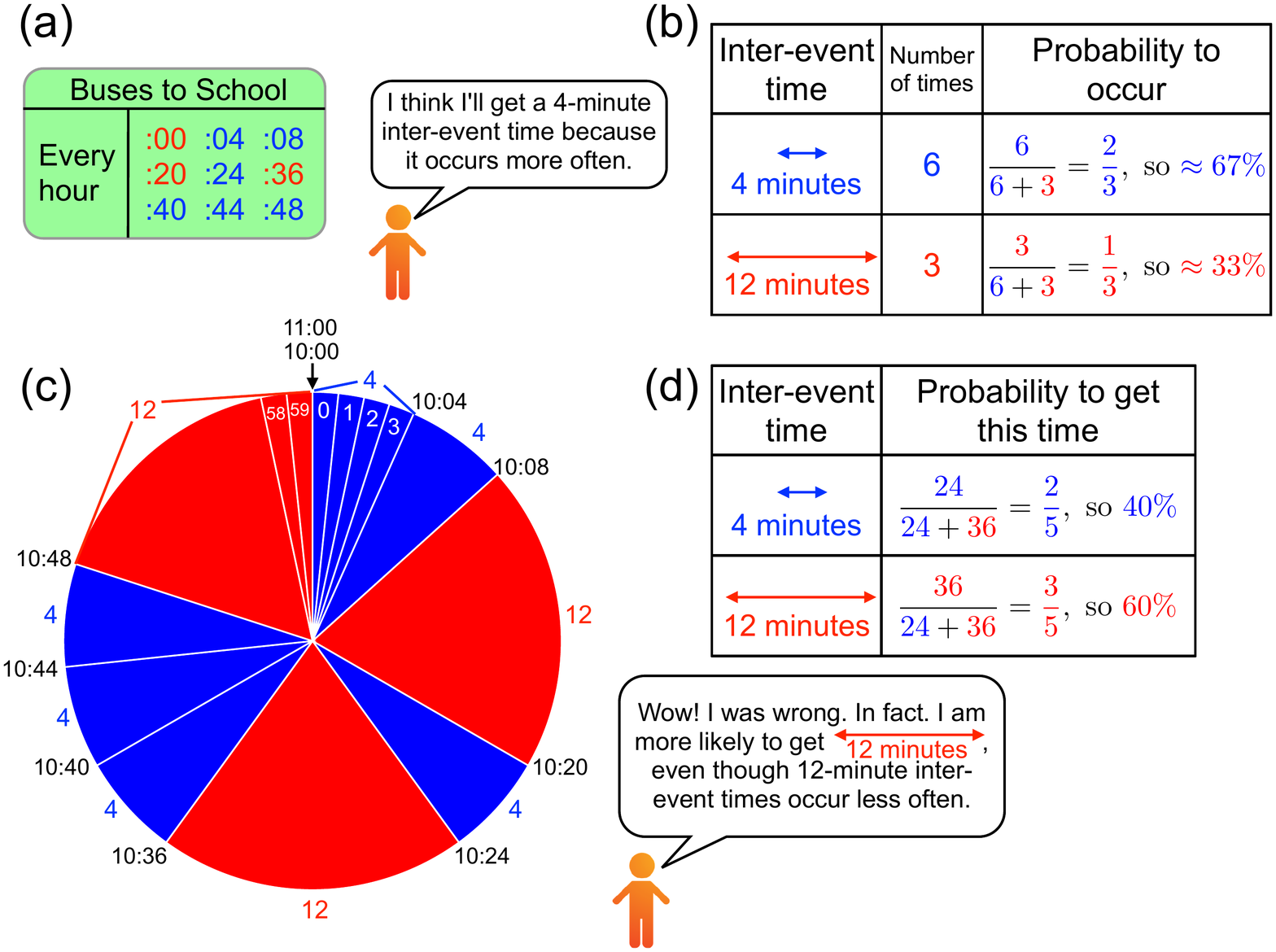}
\caption{Comparison between the naive guess of an average waiting time and the actual average waiting time in a bus schedule. (a) A simplified bus schedule.
(b) Naive guess about the probability to observe a 4-minute or 12-minute inter-event time based on the number of those inter-event times. 
Because there are 6 inter-event times that are 4 minutes long and 3
inter-event times that are 12 minutes long, there are $6 + 3 = 9$
inter-event times in total each hour. The probability that one gets a
4-minute inter-event time is $6/9$, which is about $0.67$.
(c) The organization of inter-event times in one hour. (d) The actual probability to get a 4-minute or 12-minute inter-event time. 
}
\label{fig:waiting-time paradox calculated}
\end{figure}

The bus schedule in Figure~\ref{fig:bus example}(b) looks complicated, and it has various irregular-looking numbers.
In general, when a situation looks complicated, it's useful to simplify it before applying mathematical reasoning. This is a key idea behind mathematical modeling, and it is something that professional applied mathematicians, physicists, and others do all of the time. With this in mind, let's consider the simpler bus schedule in Figure~\ref{fig:waiting-time paradox calculated}(a).
In this simpler schedule, buses are scheduled to arrive with either 4-minute or 12-minute intervals between buses. It is simpler than the original schedule in Figure~\ref{fig:bus example}(b) because there are several different intervals between bus arrivals in the original schedule. The waiting-time paradox also occurs in the simplified bus schedule, and the simplification makes it easier for us to explain --- and for you to see --- what is going on.

\subsection{You tend to get a long ``inter-event time'', which causes the paradox}

A key thing that we need to consider is the \textbf{inter-event time}, which is the time between two consecutive buses.
The schedule in Figure~\ref{fig:waiting-time paradox calculated}(a) indicates that, in each hour,
6 buses arrive immediately after an inter-event time of 4 minutes (we show this in blue) and 3 buses arrive immediately after an inter-event time of 12 minutes (we show this in red). In our scenario, recall that you have arrived at the bus stop without knowing what time it is.
You may arrive at the bus stop at 10:23, at 10:41, or at some other time. This may seem unrealistic, but even if you carry your mobile phone so that you know what time it is, if you just show up at a bus stop without checking the bus schedule, then mathematically it is the same as not knowing the time. In other words, you do not know if the next bus will come very soon or if it will take a while. In the language of probability theory, we say that you are arriving at the bus stop at a time that has been chosen \textbf{uniformly at random}.

Now that we have established that you have arrived at the bus stop at a time that we pick uniformly at random, do you think that you are more likely to have a 4-minute inter-event time or a 12-minute one? If you get a 4-minute inter-event time, you will not wait very long for the next bus. Your waiting time will be, for example, 1 minute or 3 minutes. However, if you get a 12-minute one, you may have to wait a long time, such as 10 minutes. As we indicated above, 6 inter-event times are 4-minutes long, and 3 inter-event times are 12-minutes long. 

Perhaps you will get a 4-minute inter-event time because there are more of those (there are 6 of them each hour) than there are 12-minute inter-event times (there are 3 of them each hour). The former occur with probability $6/(6+3) = 2/3$, which is about $0.67$, so this situation occurs about 67\% of the time [see
 Figure~\ref{fig:waiting-time paradox calculated}(b)]. The latter occur with probability $3/(6+3) = 1/3$, so this situation occurs about 33\% of the time.
  Unfortunately, this is wishful thinking.
Take a look at Figure~\ref{fig:waiting-time paradox calculated}(c). Arriving at the bus stop ``uniformly at random'' is like rotating a wheel with the numbers $0$,
$1$, $2$, $\ldots$, $59$ at each `slice' of the wheel (say, each in a different color) in a game at a high speed, stopping it with your finger, and looking at the 
slice that your finger is touching.
If your finger points to 33, it means that you arrive at the bus stop at 10:33. In this case, you need to wait 3 minutes for the next bus, which arrives at 10:36. The game wheel is a device to realize the notion of ``uniformly at random'' in practice. We color each minute (from $0$ to $59$) based on the 
inter-event times
in Figure~\ref{fig:waiting-time paradox calculated}(c). The coloring scheme is the same as in 
Figures~\ref{fig:waiting-time paradox calculated}(a,b), so blue corresponds to a 4-minute inter-event time and red corresponds to a 12-minute one.
There are 24 blue minutes in Figure~\ref{fig:waiting-time paradox calculated}(c), and the other 36 minutes are red. This is how 1 hour is divided into blue and red minutes. From this picture, we see that you are more likely to get a 12-minute inter-event time (with a probability of $0.60$) than a 4-minute one (with a probability of $0.40$). See Figure~\ref{fig:waiting-time paradox calculated}(d).

Only 3 of the 9 inter-event times (i.e., with probability $1/3$, which is about $0.33$) are 12 minutes long, but it is still more likely to get one of these (with probability $0.60$) than one of the 4-minute inter-event times. Why is this the case? The answer comes from a simple fact: a long inter-event time is long, and a short one is short. More specifically, a long inter-event time occupies 12 numbers of a game wheel, but a short one occupies only 4 numbers. The 3 long inter-event times together cover $12 \times 3 = 36$ of the 60 minutes.
 By contrast, the 6 short ones together cover only $4\times 6 = 24$ minutes. By the same argument, the rareness of an inter-event time in a bus schedule does not imply that it is rare to encounter that inter-event time.


\subsection{How big is the impact of the waiting-time paradox?}

The essential mathematical logic behind the waiting-time paradox is that you are more likely to get a long inter-event time. Once you do get a long inter-event time, it is likely that you will have to wait a long time for the next bus. By doing 
some more calculations with the bus schedule in Figure~\ref{fig:waiting-time paradox calculated}(a), we see that the naive guess for how long you should expect to wait (i.e., half of the average inter-event time) is 3 minutes and 20 seconds, but the actual waiting time is 4 minutes and 24 seconds on average. In the schedule in Figure~\ref{fig:bus example}(b), the naive guess for the average waiting time is 3 minutes,
 and the actual waiting time on average is 3 minutes and 40 seconds.

Figure~\ref{fig:waiting-time paradox calculated}(a) illustrates a simplistic setting. It has two values of inter-event times: 4 minutes and 12 minutes. Let's further simplify this bus schedule by supposing that all inter-event times are 6 minutes long.
Therefore, buses arrive regularly every 6 minutes on average.
In this case, there is no longer a waiting-time paradox, because
the naive guess and the actual average waiting time are both 3 minutes long. For the waiting-time paradox to occur,
we need to have a mixture of at least two inter-event times, such as 4 minutes and 12 minutes.


\section{A few applications}

\begin{figure}
\includegraphics[width=16cm]{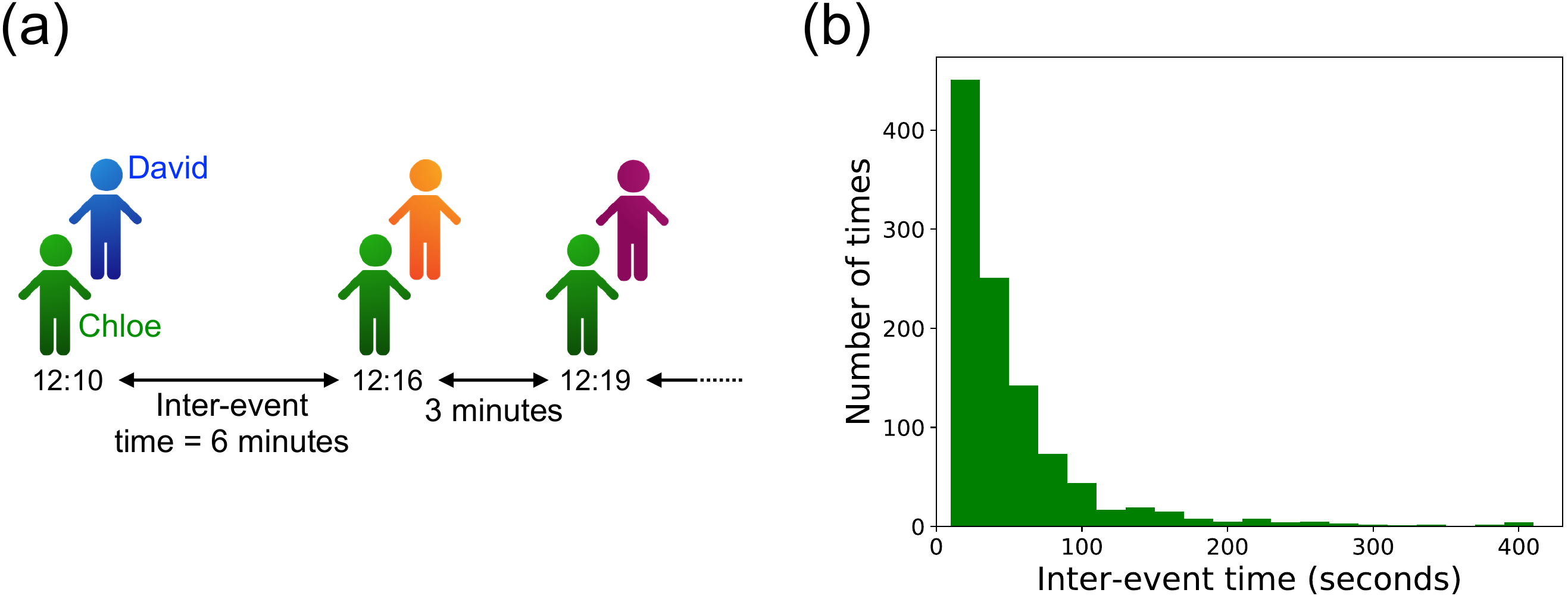}
\caption{Inter-event times in a 
school. (a) Illustration of the inter-event times for a student named Chloe. By ignoring the students with whom Chloe has talked, we conclude that two of her inter-event times are 6 minutes and 3 minutes. (b) The number of each inter-event time for a student in a school in France. This example comes from the ``Primary School'' data set in the SocioPatterns project \cite{Isella2011PlosOne}. We selected the student with the largest number of events and calculated all of their inter-event times. The bars in this picture, which together form a type of picture called a \textbf{histogram}, indicate the number of inter-event times of each length for that student. 
The data set has recordings of interactions between students 
each 20 seconds.
For example, if Chloe is talking to David 12:10:00 and again at 12:10:20, then we count the two observations as part of a single conversation (i.e., event) between them. However, if Chloe is talking to David at 12:10:00 and 12:10:40 but not at 12:10:20, we assume that there are two separate events.
}
\label{fig:primary school}
\end{figure}

The implication of the waiting-time paradox is much wider than just waiting for buses.
 Inter-event times are important for many situations. At school, consider the ``event'' of talking to one of your classmates, rather than the event of a bus arriving.
After the next time that you talk to one of your classmates (in other words, after the next event), the inter-event time is the amount of time between those conversations [see Figure~\ref{fig:primary school}(a)]. One inter-event time may be 2 minutes, and the next one may be 11 minutes. In social networks like ones that represent who talks to whom, inter-event times typically vary much more drastically than they do for bus arrivals. If buses arrive with 3-minute intervals but then an hour passes before the next bus is scheduled to arrive, people would get very upset (and many people may be late for school). 
However, in social activities, there are often large variations in inter-event times. In Figure~\ref{fig:primary school}(b), we show an example of inter-event times
for students in a school in France \cite{Isella2011PlosOne}. We list these times from left to right in increasing order in a diagram that is called a \textbf{histogram}. The height of the bar indicates the number of times that each inter-event time occurs.
Most inter-event times are short (such as 20 or 40 seconds), but a small number of them are large (such as between 200 and 400 seconds).
These kinds of large inter-event times can have a major impact on many things, like the speed of the spread of a disease in a population \cite{Karsai2011PhysRevE}.

In our example of the waiting-time paradox with buses, we saw that even if there are only 3 long inter-event times among 9 total inter-event times, we are more likely to get one of the long inter-event times than one of the short ones. This is an example of \textbf{biased sampling}.
Another famous example of biased sampling is the 
friendship paradox \cite{Feld1991AmJSociol,Strogatz2012NewYorkTimes}.
According to the friendship paradox, your friends tend to have more friends than you do. 
However, there is no reason to get upset, because this too is a purely mathematical phenomenon.
If you have 20 friends in your school, many of them are likely to be popular people.
For example, if Alice has just one friend, it is unlikely that you are Alice's only friend; it's more likely to be someone else.
By contrast, if Bob is friends with half of the students in your school, then it is very likely that you are one of Bob's friends.
Waiting for the next bus and counting the number of friends may seem to have nothing to do with each other. However, from a mathematical perspective, you are likely to have a friend like Bob for basically the same reason that you are likely to get a bus after a long inter-event time.
Suppose that there are 6 students who each have 4 friends and 3 students who each have 12 friends.
If you are friends with just one of these 10 students, then your friend is likely to be one with 12 friends, even though there are only 3 students with 12 friends among the $6+3 = 9$ students. These numbers are exactly the same as the ones that we used in Figure~\ref{fig:waiting-time paradox calculated} to demonstrate the waiting-time paradox.
This illustrates that the waiting-time paradox and the friendship paradox have the same mathematical origin, and both are a consequence of biased sampling.

As we have seen, mathematics provides a way to unify seemingly different ideas and to see when they are closely related. This is true not just with the waiting-time paradox and the friendship paradox, but also with many other things.


\section*{Glossary}

\begin{itemize}

\item \textbf{Biased sampling:} 
This occurs when one obtains (or, to phrase it more technically, one ``samples'') items, such as long inter-event times or a person with many friends, with a larger probability than other items for a systematic reason.

\item 
\textbf{Histogram:} A diagram that shows the number of items in several ranges of numbers (or in many such ranges) to compare the amounts for different ranges. For example,
suppose that one shows children with the age ranges 0--4, 5--9, 10--14, $\ldots$ from left to right. The height of the bar for the range 10--14 indicates the number of people whose ages are between 10 and 14 years old. One compares the height of this bar to those of the other bars.


\item \textbf{Inter-event time:} The length of the time between two consecutive events, such as the arrivals of two buses at a bus stop or two conversation events of a person.

\item \textbf{Probability theory:} A research field in mathematics about topics that are related to ``probability,'' which is a numerical description of how likely an event is to occur. A ``probability distribution'' is a mathematical function that gives the probabilities of occurrence of different outcomes of an event. 

\item \textbf{Uniformly at random:} A probability distribution in which each possible event is equally likely to occur.

\end{itemize}


\section*{Acknowledgements}

We are grateful to our young readers --- Nia Chiou, Taryn Chiou, Valerie E. Eng, Anthony Jin, Iris Leung, Maple Leung, Ami Masuda, and Ritvik Mukherjee --- for their many helpful comments. We also thank their parents and teachers --- Lyndie Chiou and Christina Chow --- for putting us in touch with them and soliciting their feedback. We thank Iris Leung for drawing Figure~\ref{fig:bus example}(a). We thank the SocioPatterns collaboration (see \url{http://www.sociopatterns.org}) for providing data. MAP acknowledges support from the National Science Foundation (grant number 1922952) through the Algorithms for Threat Detection (ATD) program. 



\begin{thebibliography}{1}

\bibitem{Weldman1957JOperResSoc}
Welding, P. I.
The instability of a close-interval service.
\emph{Journal of the Operational Research Society}
\textbf{8}, 133--142 (1957).

\bibitem{Masuda2020NortheastJCompSyst}
Masuda, N. \& Hiraoka, T.
Waiting-time paradox in 1922.
\emph{Northeast Journal of Complex Systems}
\textbf{2}, 1 (2020).

\bibitem{Isella2011PlosOne}
Isella, L. \emph{et~al.}
Close encounters in a pediatric ward: Measuring face-to-face proximity and mixing patterns with wearable sensors.
\emph{PLoS ONE} \textbf{6},
e17144 (2011).

\bibitem{Karsai2011PhysRevE}
Karsai, M. \emph{et~al.}
Small but slow world: How network topology and burstiness slow down spreading.
\emph{Physical Review E} \textbf{83}, 025102(R) (2011).

\bibitem{Feld1991AmJSociol}
Feld, S.~L.
Why your friends have more friends than you do.
\emph{American Journal of Sociology} \textbf{96},
1464--1477 (1991).

\bibitem{Strogatz2012NewYorkTimes}
Strogatz, S. H.
Friends you can count on.
\emph{The New York Times}
\textbf{17 September 2012} (2012).
Available at \url{https://opinionator.blogs.nytimes.com/2012/09/17/friends-you-can-count-on/}.

\end{thebibliography}


\section*{Author Biographies}


\textbf{Naoki Masuda} is an associate professor in the Department of Mathematics at the State University of New York at Buffalo. Born in Tokyo, Japan, he has lived in Tokyo, San Diego (USA), Bristol (UK), and Buffalo (USA). He used to cycle but now has jogged regularly for more than 10 years.
%
%
Cycling and running are often faster transportation methods than buses in congested cities like Tokyo and Bristol. In Great Britain, he recently learned the following saying: ``You wait forever for a bus, and then three come at once.'' On several occasions in Great Britain, Naoki saw
three buses that were bound for the same destination come in tandem. When he missed this wave of buses, he had to wait a very long time for the next bus, which is consistent with the waiting-time paradox.

\medskip

\noindent \textbf{Mason A. Porter} is a professor in the Department of Mathematics at UCLA. He was born in Los Angeles, California, and he is excited to be a professor in his hometown. In addition to studying networks and other topics in mathematics and its applications, Mason is a big fan of games of all kinds, fantasy, baseball, the 1980s, and other wonderful things. Mason used to be a professor at University of Oxford, where he did actually wear robes on occasion (like in the Harry Potter series). His main memories of the waiting-time paradox are with the delightful buses in Berkeley, California, where Mason lived when he spent a semester at the Mathematical Sciences Research Institute.

\end{document}